\documentclass[12pt]{article}
\usepackage{amsxtra,amssymb,amsthm,amsmath,latexsym}

\theoremstyle{plain}

\def\R{{\mathbb R}}

\def\oH{{\overset{\circ}{H}}}
\def\oH1{{\overset{\circ}{H}\kern-.02in{}^1}}

\def\bee{\begin{equation*}}
\def\eee{\end{equation*}}
\def\be{\begin{equation}}
\def\ee{\end{equation}}

\begin{document}

\title{ Global existence, uniqueness and estimates of the solution to the Navier-Stokes equations
}

\author{Alexander G. Ramm\\
 Department  of Mathematics, Kansas State University, \\
 Manhattan, KS 66506, USA\\
ramm@math.ksu.edu\\
}

\date{}
\maketitle\thispagestyle{empty}

\begin{abstract}

The Navier-Stokes (NS) problem consists of finding a vector-function $v$ from the Navier-Stokes equations.
The solution $v$ to NS problem is defined in this paper as the solution to an integral equation. The kernel $G$ of this equation
solves a linear problem which is obtained from the NS problem by dropping the nonlinear term $(v \cdot \nabla)v$.
The kernel $G$ is found in closed form. Uniqueness of the solution to the integral equation is proved
in a class of solutions $v$ with finite norm $N_1(v)=\sup_{\xi\in \R^3, t\in [0, T]}(1+|\xi|)(|v|+|\nabla v|)\le c   (*)$,
where $T>0$ and $C>0$ are arbitrary large fixed constants. In the same class of solutions existence of the solution is proved
under some assumption.
 Estimate of the energy of the solution is given.

\end{abstract}

\section{Introduction}\label{S:1}
There is a large literature on the Navier-Stokes (NS) problem, ( see \cite{L}, \cite{T} and references therein,
including the papers by Leray, Hopf, Lions, Prodi, Kato and others). The global existence and uniqueness of a solution
was not proved.
 Various notions of the solution were used, see \cite{L}, \cite{LR}, \cite{R646}. In this paper a new notion of the solution is given:
the solution to NS problem is defined as a solution to some integral equation. It is proved that this solution is unique, it exists
 in a certain class of functions globally, i.e.,  for all $t\ge 0$, it has finite energy, and an estimate of this solution is given. 
 
  The NS problem consists of solving the equations
\be\label{e1} v'+(v, \nabla)v=-\nabla p +\nu\Delta v +f, \quad x\in \R^3,\,\, t\ge 0,\quad \nabla \cdot v=0, v(x,0)=v_0(x).
\ee
Vector-functions $v=v(x,t)$, $f=f(x,t)$ and the scalar function $p=p(x,t)$ decay as $|x|\to \infty$ uniformly with respect to
$t\in [0, T]$, where $T>0$ is an arbitrary large fixed number,  $v'=v_t$,
$\nu=const>0$, $v_0$ is given, $ \nabla \cdot v_0=0$, the velocity $v$ and the pressure $p$ are unknowns,  $v_0$ and $f$ are known.
Equations \eqref{e1}  describe viscous incompressible fluid with density $\rho=1$.
We assume that $|f|+|\nabla f|$ and $|v_0|+|\nabla v_0|$ decay as $O(\frac 1 {(1+|x|)^{a}})$ as $|x|\to \infty$, where $a>3$
and for $f$ the decay with respect to $x$ holds for any $t\ge 0$,  uniformly with respect to $t\in [0, T]$,  where  $T\ge 0$ is an arbitrary large fixed number.

 Our approach consists of three steps.
First, we construct tensor $G_{jm}(x,t)$, equivalent to Oseen's thensor in \cite{LR}, p.62, solving the problem:
\be\label{e2} G'-\nu \Delta G=\delta(x)\delta(t)\delta_{jm}- \nabla p_m,\quad \nabla \cdot G=0; \quad G=0, \quad t<0. \ee
Here $\delta(x)$ is the delta-function, $\delta_{jm}$ is the Kronecker delta, $p_m=p(x,t)e_m$, $p$ is the scalar,
$p(x,t)=0, t<0$; $\nabla p_m=\frac {\partial p}{\partial x_j}e_je_m$,
$\{e_j\}|_{j=1}^3$ is the orthonormal basis of $\R^3$, $e_je_m$ is tensor. Our construction method differs from the one in \cite{LR}.

Secondly, we prove that solving \eqref{e1}
is equivalent to solving the following integral equation, cf \cite{LR}, p.62:
\be\label{e3} v(x,t)=\int_0^tds \int_{\R^3} G(x-y,t-s)[f(y,s)-(v, \nabla)v]dy +\int_{\R^3} g(x-y,t)v_0(y)dy,
\ee
where
\be\label{e3'}
 g(x,t)=\frac{e^{-\frac {|x|^2}{4\nu t}}}{(4\nu t \pi)^{3/2}}, 
 \quad g(x,0)=\delta(x), \quad g'-\nu \Delta g=\delta(x) \delta(t);\quad \int_{\R^3}g(x,t)dx=1.
\ee
Thirdly, we prove by a new method,  using
the new assumption \eqref{e19a}, see below, that equation \eqref{e3} has a solution  in the space $X$ of functions with finite norm
$N_1(v):=\sup_{x\in \R^3, t\in [0,T]}\{(|v(x,t)|+|\nabla v(x,t)|)(1+|x|)\}$ and this solution is unique in $X$. In Lemma 1 a new a priori estimate is obtained.

{\bf Theorem 1.} {\em Problem \eqref{e1} has a unique solution in $X$. The solution to \eqref{e1} exists in $X$
provided that $\sup_{n>1} N_1(v_n)< \infty$, where $v_n$ is defined in \eqref{e19}.
This solution has finite energy $E(t)=\int_{\R^3}|v(x,t)|^2dy$ for every $t\ge 0$.}

\section{Construction of $G$}\label{S:2}

Let $G=\int_{\R^3}H(\xi,t)e^{i\xi \cdot x}d\xi$. Taking Fourier transform of \eqref{e2} yields
\be\label{e4}
H'+\nu \xi^2 H=\frac {\delta(t)\delta_{jm}}{(2\pi)^3} - i\xi P_m(\xi,t), \quad \xi H=0; \quad H(\xi,t)=0,\, t<0.
\ee
The $H=H_{jm}$ is a tensor, $1\le j,m\le 3$, $P_m(\xi,t):=(2\pi)^{-3}\int_{\R^3}p_m(x,t)e^{-i\xi \cdot x}dx,$
$\xi=\xi_je_j$, summation is understood here and below over the repeated indices, $1\le j \le 3$, $\xi H:=\xi_j H_{jm}$. Multiply \eqref{e4} by $\xi$ from the left, use the equation $\nabla \cdot G=0$ which implies $\xi H=0$, $\xi H'=0$, and get
\be\label{e5}
(2\pi)^{-3}\delta(t) \xi_m=i\xi^2 P_m, \quad P_m:=P_m(\xi, t)=-i \frac {\delta(t)\xi_m} {(2\pi)^{3} \xi^2} ,
\quad \xi^2:=\xi_j \xi_j.
\ee
From \eqref{e4}-\eqref{e5} one obtains:
\be\label{e6} H_{jm}'+\nu \xi^2 H_{jm}=(2\pi)^{-3}\delta (t) \Big(\delta_{jm}- \frac{\xi_j \xi_m} {\xi^2}\Big).
\ee
Thus,
\be\label{e7} H=H_{jm}(\xi,t)=(2\pi)^{-3} \Big(\delta_{jm}- \frac{\xi_j \xi_m} {\xi^2}\Big) e^{-\nu \xi^2 t}; \quad H=0, \,\, t<0.
\ee
\be\label{e8}
G(x,t)= (2\pi)^{-3}\int_{\R^3} e^{i\xi \cdot x} \Big(\delta_{jm}- \frac{\xi_j \xi_m} {\xi^2}\Big) e^{-\nu \xi^2 t} d\xi
:=I_1+I_2.
\ee
The integrals $I_1,I_2$ are calculated in the Appendix:
\be\label{e9}
I_1= \delta_{jm}g(x,t);\quad g(x,t):=\frac{e^{-\frac {|x|^2}{4\nu t}}}{(4\nu t \pi)^{3/2}},  \quad t>0; \quad g=0, \,\, t<0,
\ee
\be\label{e10}
I_2=\frac {\partial^2}{\partial x_j \partial x_m} \Big( \frac 1 {|x|2\pi^{3/2}} \int_0^{|x|/(4\nu t)^{1/2}} e^{-s^2}ds \Big), \quad t>0;\quad I_2=0, \,t<0.
\ee
Therefore,
\be\label{e11}
G(x,t)=\delta_{jm} g(x,t) +\frac {1} {2\pi^{3/2}}\frac {\partial^2}{\partial x_j \partial x_m}
\Big( \frac {1} {|x|} \int_0^{|x|/(4\nu t)^{1/2}} e^{-s^2}ds \Big), \,\, t>0.
\ee
\section{Solution to integral equation for $v$ satisfies NS equations }\label{S:3}
Apply the operator $L:=\frac {\partial} {\partial t} -\nu \Delta $ to the left side of \eqref{e3} and use \eqref{e2} to get
\be\label{e12}
\begin{split}
Lv=\int_0^tds\int_{\R^3}\Big[\delta(x-y)\delta(t-s)\delta_{jm}+ \partial_j \Big(\partial_m\big(\frac {1} {2\pi^{3/2}|x-y|} \int_0^{|x-y|/(4\nu t)^{1/2}} e^{-s^2}ds\big)\Big)\Big]\times \\ \Big( f(y,s)-(v(y,s), \nabla)v(y,s)\Big)dy
 =f(x,t)-(v(x,t), \nabla)v(x,t) -\nabla p_m(x,t),
\end{split}
\ee
where
\be\label{e12a}
p_m(x,t)=-\int_0^tds\int_{\R^3}\partial_m\Big(\frac {1} {2\pi^{3/2}|x-y|} \int_0^{|x-y|/(4\nu t)^{1/2}} e^{-s^2}ds\Big)
 [ f(y,s)-(v(y,s), \nabla)v(y,s)]dy,
 \ee
 and $v(x,0)=v_0(x)$ because $g(x-y,0)=\delta(x-y)$. Using the formula $\nabla \cdot G=0$, the relation
 $\int_{\R^3}g(x-y,t)v_0(y)dy=\int_{\R^3}g(z,t)v_0(z+x)dz$ and the formula $\nabla \cdot v_0=0$,
  one checks that $\nabla \cdot v=0$. Thus, a solution to \eqref{e3} solves \eqref{e1}.

\section{Proof of the uniqueness of the solution to $\eqref{e3}$ in the space of functions
with finite norm  $N_1(v)$}\label{S:4}
Let $X$ be the space  $C^1(\R^3; 1+|x|)$ of vector-functions with the norm
$$N_1(v):=\sup_{x\in \R^3, t\in [0,T]} [(|v(x,t)|+|\nabla v(x,t)|)(1+|x|)],$$
where  $T> 0$ is an arbitrary large fixed number,
 $\nabla$ stands for any first order derivative, and let $N_0(v)=\sup_{x\in \R^3} (|v(x,t)|+|\nabla v(x,t)|)$,
$N_0(v)\le N_1(v)$, $|v|:=(\sum_{j=1}^3 |v_j|^2)^{1/2}$.
Assume that there are two solutions $v_1,v_2$ to equation \eqref{e3} with finite norms $N_1(v)$ and let $w=v_1-v_2$.
 One has
\be\label{e13}
\begin{split}
w(x,t)=\int_0^tds\int_{\R^3} G(x-y,t-s)[(v_2, \nabla)v_2-(v_1, \nabla)v_1]dy&=\\
\int_0^tds\int_{\R^3} G(x-y,t-s)[(w, \nabla)v_2+(v_1,\nabla)w]dy.
\end{split}
\ee
Therefore
\be\label{e14} N_0(w)\le \int_0^tdsN_0\Big(\int_{\R^3} |G(x-y,t-s)|(1+|y|)^{-1}dy\Big)N_0(w) (N_1(v_2)+N_1(v_1))\le
c \int_0^t\frac{N_0(w)}{ (t-s)^{ 1/ 2}}ds.
\ee
 Since we are proving uniqueness in the set of functions with finite norm  $N_1(v)$,
one has $N_1(v_2)+N_1(v_1)\le c$, where $c>0$
stands for various estimation constants. It is checked in the Appendix that
\be\label{e14'}
N_0\Big(\int_{\R^3} |G(x-y,t)|(1+|y|)^{-1}dy\Big)\le ct^{-\frac 1 2},
 \ee
 so inequality \eqref{e14} holds.
 From \eqref{e14} by the standard argument one derives that $N_0(w)=0$. Thus, $v_1=v_2$.
 Uniqueness of the solution to \eqref{e3} is proved in $X$. \hfill$\Box$
\section{Proof of the existence of the solution  to $\eqref{e3}$ in the set of functions
with finite norm  $N_1(v)$}\label{S:5}
Rewrite \eqref{e3} as
\be\label{e17} v(x,t)=-\int_0^tds\int_{\R^3} G(x-y,t-s)(v, \nabla)v dy + F(x,t),
\ee
where $F$ is known:
\be\label{e18} F(x,t):=\int_0^tds\int_{\R^3} G(x-y,t-s)f(y,s)dy +\int_{\R^3} g(x-y,t)v_0(y)dy.
\ee
Equation \eqref{e17} is of Volterra type, nonlinear,  solvable, as we wish to prove, by iterations:
 \be\label{e19} v_{n+1}(x,t)=-\int_0^t B(v_n, \nabla)v_n ds+ F, \quad v_1=F;\qquad Bv:= \int_{\R^3} G(x-y,t-s)v(y,s)dy.
 \ee
Using \eqref{e14'}, the formula $\int_0^t \frac {ds}{s^{1/2}}=2t^{1/2}$, and assuming that
\be\label{e19a}
\sup_{n\ge 1} N_1(v_n)\le c,
\ee
 one gets:
\be\label{e20}
\begin{split}
 N_0(v_{n+1}-v_n)\le N_0\Big(\int_0^tds\int_{\R^3} \frac{|G(x-y,t-s)|}{1+|y|}dy\Big)[N_0(v_n-v_{n-1})(N_1( v_n)+N_1(v_{n-1})] &\le\\
 ct^{1/2}N_0(v_n-v_{n-1}).
\end{split}
\ee
 Therefore,
 \be\label{e21}
 N_0(v_{n+1}-v_n)\le ct^{1/2}N_0(v_n-v_{n-1}).
 \ee
 If $\tau $ is chosen so that $ c\tau^{1/2}<1$, then $B$ is a contraction map on the bounded set $N_0(v)\le R$ for sufficiently large $R$ and $t\in [0,\tau]$. Thus, $v$ is uniquely determined by iterations for $x\in \R^3$ and $t\in [0,\tau]$. If $t>\tau$ rewrite \eqref{e17} as
  \begin{align*} v=F- \int_0^{\tau}ds\int_{\R^3} G(x-y,t-s)(v, \nabla)v dy-\int^t_{\tau}ds\int_{\R^3} G(x-y,t-s)(v, \nabla)vdy &:=\\ F_1-\int^t_{\tau}ds\int_{\R^3} G(x-y,t-s)(v, \nabla)v dy,
  \end{align*}
  where $F_1$ is a known function since $v$ is known for  $t\in [0,\tau]$. To this equation one applies the contraction mapping principle and get $v$ on the interval $[0, 2\tau]$. Continue in this fashion and construct $v$ for any $t\ge 0$.
 Process \eqref{e19} converges to a solution to  equation \eqref{e3}.  \hfill$\Box$

 Let us make a remark concerning the mapping done by the operator in \eqref{e17}. In \cite{MP}, p. 234, sufficient conditions are given for a singular integral operator to map a class of functions with a known power rate of decay at infinity into a class of functions with a suitable rate of decay at infinity.
It follows from \cite{MP}, Theorem 5.1 on p.234, that if $|f(x)|+|\nabla f|\le c(1+|x|)^{-a}$, $a>3$,  then the part of the operator $G$ in \eqref{e17},  responsible for the lesser decay of the iterations $v_n$ at infinity,  acts as a weakly singular operator similar to the operator $Q$,  $Qf:=\int_{\R^3} \frac {b(|x-y|)|f(y)|}{|x-y|^3}dy \le c \frac{1}{(1+|x|)^3}$ for large $|x|$. This part of $G$ yields the decay $O(\frac{1}{(1+|x|)^3})$ at infinity. This decay, in general,  cannot be improved.
 The first iteration $v_1$ yields a function decaying
with its first derivatives as $O(\frac{1}{(1+|x|)^3})$. The second iteration contains a function whose behavior for large $|x|$
is determined by the decay of the functions $F$. Since $(v_1\cdot \nabla)v_1=O((1+|x|)^{-7})$, the decay of $v_2$ is again $O((1+|x|)^{-3})$
because $v_1=O((1+|x|)^{-3})$ and $|\nabla v_1|=O((1+|x|)^{-4})$. Thus, for $|v_n|+|\nabla v_n|$ one gets the decay of the order
$ O((1+|x|)^{-3})$ as $|x|\to \infty$ for every fixed $t\ge 0$, provided that $f$ and $v_0$ together with their first derivatives decay  not slower than $O((1+|x|)^{-a})$, $a>3$.

 \section{Energy of the solution}\label{S:6}
 In this Section we prove that the solution has finite energy in a suitable sense and give an estimate of the solution as $t\to \infty$. Let us define the energy by the integral $E(t)=\int_{\R^3}|v(x,t)|^2dy$. If one multiplies \eqref{e1} by $v$ and
 integrate over $\R^3$, then one gets a known conservation law (see \cite{L}):
 \be\label{e22}
\frac 1 2 \frac{\partial}{\partial t}  \int_{\R^3}|v(x,t)|^2dx +\nu \int_{\R^3}\nabla v_j(x,t)\cdot\nabla v_j(x,t)dx= \int_{\R^3}f\cdot v dx.
 \ee
 Integrating \eqref{e22} with respect to $t$ over any finite interval $[0,t]$, $0\le t \le T$, and denoting
 $\mathcal{N}(v):= (\int_{\R^3}|v(x,t)|^2dy)^{1/2}=\sqrt{E(t)}$, one gets
 \be\label{e23}
 2\nu \int_0^t \|\nabla v\|^2dt+E(t)\le E(0) +2\int_0^t ds | \int_{\R^3}f\cdot v dx|\le E(0)+2\int_0^T \mathcal{N}(f)\mathcal{N}(v)ds.
 \ee
 Denote $E_T:=\sup_{t\in [0,T]}E(t)$. Maximizing inequality \eqref{e23} with respect to $t\in [0,T]$, and using the
 elementary inequality $2ab\le \epsilon a^2+\frac 1{\epsilon} b^2$, $a,b, \epsilon>0$, one derives from \eqref{e23} the following inequality
 \be\label{e24}
  (1-\epsilon)E_T\le  E(0)+ \frac 1{\epsilon} (\int_0^T \mathcal{N}(f)dt)^2.
  \ee
  Let $\epsilon=\frac 1 2$. Then \eqref{e23}-\eqref{e24}  allow one to estimate $E_T$ and $ \nu \int_0^t \|\nabla v\|^2dt$ through $f$ and $v_0$.
  In particular, we have proved the following theorem

  {\bf Theorem 2.} {\em Assume that $\int_0^T \mathcal{N}(f)dt<\infty$ for all $T>0$ and $E(0)<\infty$. Then $\sup_{t\in [0,T]}E(t)<\infty$
  and $\int_0^T \|\nabla v(x,t)\|^2dt<\infty$ for all $T>0$.}

 {\bf Remark 1.} It is proved in \cite{R621} that in a bounded domain $D\subset \R^3$
   the solution to a boundary problem for equations \eqref{e1} in a bounded domain with the Dirichlet boundary condition for large $t$ decays exponentially provided that $\int_0^\infty e^{bt}\mathcal{N}(f)dt<\infty$ for some $b=const>0$.

 \section{Appendix}\label{S:7}
 1. {\em Integral $I_1$.}  One has $I_1=I_{11}I_{12}I_{13}$, where
  $I_{1p}:=\frac 1 {2\pi}\int_{-\infty}^{\infty}e^{ix_p\xi_p -\nu \xi^2t}d\xi_p.$ No summation in $ x_p\xi_p$ is understood here.
   One has
   $$\frac 1 \pi\int_0^\infty \cos(x_p\xi_p)e^{-\nu \xi^2t}d\xi_p=\frac 1 \pi \frac {\pi^{1/2}}{2(\nu t)^{1/2}} e^{ -\frac {x_p^2}{4\nu t}}=\frac {e^{-\frac{x_p^2}{4\nu t }}}{(4\nu t \pi)^{1/2}}.
   $$
   After a multiplication (with $p=1,2,3$)  this yields formula \eqref{e9}.

{\em Integral $I_2$.} One has
  \be\label{e44} I_2= \frac {\partial^2}{\partial x_j \partial x_m}\frac 1 {(2\pi)^3}\int_{\R^3} e^{i\xi \cdot x} \frac{e^{-\nu t \xi^2}}{\xi^2}d\xi=\frac {\partial^2}{\partial x_j \partial x_m}\Big(\frac 1 {(2\pi)^2}\int_0^\infty dr e^
 {-\nu t r^2}\int_{-1}^1 e^{ir|x|u}du\Big),
 \ee
 so
 \be\label{e45}
 I_2= \frac {\partial^2}{\partial x_j \partial x_m}\Big(\frac 1 {|x|} \frac 1 {2\pi^{3/2} }
 \int_0^{\frac {|x|}{(4\nu t)^{1/2}}} e^{-s^2}ds \Big).
 \ee
 Here we have used formula (2.4.21) from \cite{BE}:
 \be\label{e46}
 \int_0^\infty \frac{\sin (sy)}{s}e^{-as^2}ds=\frac {\pi}2 Erf(\frac y {2a^{1/2}}); \quad Erf(y):=\frac 2 {\pi^{1/2}}\int_0^y e^{-s^2}ds.
\ee
 2. One has $\int_{\R^3}|\nabla g(x,t)|dx\le (4\nu t \pi)^{-\frac 3 2}\int_{\R^3}e^{-\frac{|x|^2}{4\nu t}}\frac{2|x|}{4\nu t}dx= \frac 2 {(\nu t \pi)^{\frac 1 2}}.$ This proves the estimate\\
 $\int_{\R^3} |G(x-y,t)(1+|y|)^{-1}dy\le ct^{-\frac 1 2}$ for
 the first term of $G$ in \eqref{e11}, namely for $g$. The second term of $G$ is   $J:=\frac {\partial^2}{\partial x_j \partial x_m}
\Big( \frac{1}{|x|} \int_0^{|x|/(4\nu t)^{1/2}} e^{-s^2}ds \Big)$ up to a factor $\frac 1 {2\pi^{3/2}}$.

One checks by direct differentiation that
\be\label{e47}
\begin{split}
J=e^{-\frac{|x|^2}{4\nu t}}\Big( \frac{\delta_{jm}}{|x|^2 (4\nu t)^{1/2}} -\frac{3x_jx_m}{|x|^4 (4\nu t)^{1/2}} -\frac{2x_jx_m}{|x|^2 (4\nu t)^{3/2}}\Big)+\int_0^{|x|/(4\nu t)^{1/2}} e^{-s^2}ds \Big(\frac{3x_jx_m}{|x|^5} -\frac{\delta_{jm}}{|x|^3}\Big).
\end{split}
\ee
Let $b(x):=b(\frac {|x|}{(4\nu t)^{1/2}}):=\int_0^{|x|/(4\nu t)^{1/2}} e^{-s^2}ds$. Note that $b(s)=O(s)$ as $s\to 0$ and $b(s)\to \frac{\pi^{1/2}}{2}$ as $s\to \infty$. It is sufficient to estimate the term with the strongest singularity, namely, \newline
 $J_1:=\int_{\R^3} \frac{b(|x-y|)}{|x-y|^3 (1+|y|)}dy$.
We prove that $\max_{x\in \R^3}J_1\le c$. The other two terms in the second term of $J$ in \eqref{e47} are estimated similarly.
One has $J_1=\int_{\R^3} \frac{b(|z|)}{|z|^3 (1+|x+z|)}dy=\int_{|z|\le 1}+\int_{|z|\ge 1}:=J_{11}+J_{12}$, where $J_{11}\le c(1+|x|)^{-1},\,\, x\in \R^3$, because $b(|z|)|z|^{-3}\le c|z|^{-2}$ for $|z|\le 1$. Let us estimate $J_{12}$. Note that $1+|y|^a\le (1+|y|)^a\le c(1+|y|^a)$ for $a\ge 1$. We prove that $J_{12}\le  c\frac {\ln (3+|x|)}{2+|x|}$  for all $x\in \R^3$.
Use the spherical coordinates  $|z|=r$, $|x|:=q$, the angle
 between $x$ and $z$ denote by $\theta$, $\cos \theta=u$, and get for $a=1$:
 \be\label{e48}
 \begin{split}
J_{12}=2\pi\int_1^\infty \frac {dr}{r}\int_{-1}^1\frac 1{1+(r^2+q^2-2rqu)^{1/2}}du\le c\int_1^\infty \frac {dr}{r}\int_{-1}^1\frac 1{(1+r^2+q^2-2rqu)^{1/2}}du\\ \le cq^{-1} \int_1^\infty \frac {dr}{r^2}\frac {(r+q)^2-(r-q)^2}{(1+(r+q)^2)^{1/2}+(1+(r-q)^2)^{1/2}}\le
c  \int_1^\infty \frac {dr}{r(2+r+q)}\le  c \frac {\ln (3+q)}{2+q}.
\end{split}
\ee
{\bf Remark 2.}  In \cite{L} it is shown that the smoothness properties of the solution $v$
are improved when the smoothness properties of $f$ and $v_0$  are improved.
This also follows from equation \eqref{e3} by the properties of weakly singular integrals and
 of the kernel $g(x,t)$.

Let us prove {\em a new a priori estimate}.

{\bf Lemma 1.} {\em If the assumptions of Theorem 2 hold, then $(1+\xi^2)|\tilde{v}(\xi,t)|^2\le c+ ct$ for all $t\ge 0$, where $c>0$ does not depend on $\xi, t$.}

{\em Proof.} The Fourier transform of \eqref{e17} yields
 \be\label{e51}
\tilde{v}(\xi,t)=\tilde{F}(\xi,t)-\int_0^tds H(\xi,t-s) \widetilde{(v,\nabla)v},\quad v(x,t)=\int_{\R^3}e^{i\xi \cdot x}\tilde{v}(\xi,t)d\xi,
\ee
where $H$ is defined in \eqref{e7},
$$ |\widetilde{(v,\nabla)v}|=|\tilde{v_j}\ast i\xi_j\tilde{v}_m|\leq \mathcal{N}(\tilde {v}) \mathcal{N}(|\xi|\tilde {v}),$$
over the repeated indices $j$ one sums up, $\mathcal{N}(\tilde{v}):=\|\tilde{v}\|_{L^2(\R^3)}$, and $\tilde{v}\ast \psi$ denotes the convolution
of two functions,  so
\be\label{e52}
|\widetilde{(v,\nabla)v}|\le \mathcal{N}(\tilde {v})\mathcal{N}(|\xi| \tilde {v})\le c \mathcal{N}(|\xi|\tilde {v}),
\ee
where the known estimate $\mathcal{N}(\tilde {v})\le c$ was used. One has
\be\label{e53} |H|\le ce^{-\nu t \xi^2}.
\ee
It follows from \eqref{e22} and \eqref{e23} that
\be\label{e54}
(2\pi)^3\int_0^tds\int_{\R^3}|\xi|^2 |\tilde{v}|^2 d\xi= \int_0^tds\int_{\R^3}|\nabla v|^2dx\le c, \quad \forall t\ge 0,
\ee
provided that assumptions of Theorem 2 hold.  It follows from \eqref{e51} that
\be\label{e55}
|\tilde{v}|\le |\tilde{F}|+c\int_0^t dse^{-\nu (t-s) \xi^2} \mathcal{N}(|\xi|\tilde{v})\le |\tilde{F}|+
c\Big(\int_0^te^{-2\nu (t-s) \xi^2}ds\Big)^{1/2} \Big(\int_0^t ds  \mathcal{N}^2(\nabla v)\Big)^{1/2}.
\ee
Thus,
\be\label{e56}
|\xi^2||\tilde{v}|^2\le 2|\xi^2||\tilde{F}|^2 +2c^2 |\xi^2|\frac {1- e^{-2\nu t \xi^2}}{2\nu \xi^2}\le c+ct,
\ee
provided that $|\xi^2||\tilde{F}|^2\le c$. This inequality holds if $f$ and $v_0$ are smooth and rapidly decaying.
Lemma 1 is proved. \hfill$\Box$


\end{document}